
\input amstex
\documentstyle{amsppt}
\nologo
\define\longdashrightarrow{{\,\text{- - - -}{\scriptscriptstyle{}^{{}_{>}}}}}
\define\longdashleftarrow{{{\scriptscriptstyle{}^{{}_{<}}}\text{- - - -}\,}}
\define\CDdashright#1#2{&\,\mathop{\longdashrightarrow}\limits^{#1}_{#2}\,&}
\define\CDdashleft#1#2{&\,\mathop{\longdashleftarrow}\limits^{#1}_{#2}\,&}
\define\CDsubset{&\quad\mathop{\subset}\quad&}
\def\Z{{\Bbb Z}}

\def\P{{\Bbb P}}

\def\C{{\Bbb C}}

\def\Til#1{\widetilde{#1}}

\def\PGL{\text{PGL}}

\def\oc{{O_C}}
\def\ocbar{\overline{O_C}}
\def\im{\text{im\,}}
\def\rk{\text{rk\,}}
\def\cP{\check{\P}}
\def\phy{\varphi}
\def\vq{V^{(q)}}
\def\bq{B^{(q)}}
\def\cq{c^{(q)}}
\def\Bq{\Cal B^{(q)}}
\def\Eq{E^{(q)}}

\def\tablerule{\noalign{\hrule}}
\def\statable#1.{$$\text{
\vbox{\tabskip=0pt \offinterlineskip
\halign to305pt{\strut##& \vrule##\tabskip=1em plus2em&
  \hfil##\hfil& \vrule\thinspace\vrule##& \hfil##\hfil& \vrule##&
  \hfil##\hfil& \vrule##& \hfil##\hfil& \vrule##&
  \hfil##\hfil& \vrule##& \hfil##\hfil& \vrule##&
  \hfil##\hfil& \vrule##\tabskip=0pt\cr\tablerule
#1}}}$$}
\def\protable#1.{$$\text{
\vbox{\tabskip=0pt \offinterlineskip
\halign to300pt{\strut##& ##\tabskip=1em plus2em&
  \hfil##\hfil& \vrule\thinspace\vrule##&
  \hfil##\hfil& \vrule##&
  \hfil##\hfil& ##\tabskip=0pt\cr
#1}}}$$}

\pagewidth{12.45 cm}\pageheight{20.35 cm}
\topmatter
\title Linear orbits of smooth plane curves
\endtitle
\author Paolo Aluffi \\ Carel Faber\endauthor
\address {Mathematics Department, Oklahoma State University, Stillwater OK
74078, USA}
\endaddress
\address {Max-Planck-Institut f\"ur Mathematik,
Gottfried-Claren-Stra\ss e 26, W-5300 Bonn 3, Germany}
\endaddress
\endtopmatter
\document

\heading \S 1. Introduction\endheading

Consider a general codimension-8 linear subspace of the $\P^{14}$
parametrizing plane quartic curves. There is a generically
finite dominant
rational map from this $\P^6$ to the moduli space of curves of genus 3;
what is the degree of this map?

To approach this kind of questions, we embark in this paper on the
study of the natural action of $\PGL(3)$, the group of automorphisms
of $\P^2$, on the projective space $\P^N$ parametrizing plane curves
of degree $d$ (thus $N=d(d+3)/2$). We are concerned here with the
orbits of (points corresponding to) smooth plane curves $C$ of degree
$d\ge3$.  These orbits $\oc$ are 8@-dimensional quasi-projective (in
fact, affine) varieties. Their closures $\ocbar$ (in $\P^N$) are
8@-dimensional projective varieties, and one easily understands that
the answer to the above question is nothing but the {\sl degree\/} of
$\ocbar$ for the general plane quartic curve $C$. In \S2 we explicitly
construct a resolution of these varieties, which we use in \S3 to
compute their degree (for every smooth plane curve $C$ of degree $d\ge
3$). In fact the construction gives more naturally the so-called {\sl
predegree\/} of the orbit closure $\ocbar$: that is, the degree of
$\ocbar$ multiplied by the order of the $\PGL(3)$-stabilizer of $C$.
It turns out that, for a smooth curve $C$ of degree $d$, the predegree
depends only on $d$ and the nature of the flexes of $C$. As is
illustrated in \S3.6, this has nice consequences related to the
automorphism groups of smooth plane curves.

We now describe the contents of this paper more precisely. Associated
with each plane curve $C$ is a natural map $\PGL(3)\to \P^N$ with image
$\oc$; we view this as a rational map $\phi_C$ from the $\P^8$ of
$3\times3$ matrices to $\P^N$. In \S2 our object is to resolve this
map by a sequence of blow-ups over $\P^8$, in fact constructing a
non-singular compactification of $\PGL(3)$ that dominates $\ocbar$.
For a smooth $C$, the base locus of $\phi_C$ is a subvariety of $\P^8$
isomorphic to $\P^2\times C$, thus smooth; after blowing up this
support, we find that the support of the base locus of the induced
rational map from the blow-up to $\P^N$ is again smooth: and we choose
it as the center of a second blow-up. Just continuing this process
(which requires a fair amount of bookkeeping) gives a good resolution
of the map. We find that the number of blow-ups needed equals the
maximum order of contact of $C$ with a line: so, for example, three
blow-ups suffice for the general curve.

Having resolved the map $\phi_C$, we compute in \S3 the predegree of
the closure of the orbit of $C$ as the 8@-fold self-intersection of
the pull-back of the class of a `point-condition', i.e., a hyperplane
in $\P^N$ parametrizing the curves passing through a given point of
$\P^2$. The main tool is an intersection formula for blow-ups from
\cite{Aluffi1}; to apply this formula, we extract from the geometry of
the blow-ups detailed in \S 2 the relevant intersection-theoretic
information, and particularly the normal bundles and intersection
rings of the centers of the various blow-ups. This leads to explicit
formulas for the predegree of $\ocbar$ in terms of the degree of $C$
and of four numbers encoding the number and type of the flexes of $C$.
For example, the answer to the question posed in the beginning (that
is, the degree for a general quartic) is \oldnos{14,280} ($d=4$ in the
Corollary in \S 3.5).

Besides the applications to automorphism groups of plane curves
already alluded to, and more examples given in \S3.6 (e.g., we compute
the degree of the trisecant variety to the $d$-uple embedding of
$\P^2$ in $\P^N$), the computation of the degree of $\ocbar$ also has
some enumerative significance: it gives the number of translates of
$C$ that pass through 8 points in general position. On a more global
level, the degree of the orbit closure of a general plane curve of
degree $d$ equals the degree of the natural map from a general
codimension@-8 linear subspace of $\P^N$ to the moduli space of smooth
plane curves of degree $d$.  A study of the boundaries of orbits and
of 1@-dimensional families of orbits reveals where this map is proper:
these matters will be treated in a sequel to this paper. Also, we hope
to be able to unify some scattered results we have concerning orbits
of {\sl singular\/} curves; and we plan to study the singular locus of
orbit closures.

Excellent practice to become familiar with the techniques of the
paper is to apply them to the easier case of the action of the group
$\PGL(2)$ on the spaces $\P^d$ parametrizing $d$-tuples of points on a
line. Only one blow-up of the $\P^3$ of $2\times 2$ matrices is needed
in this case, and one finds the following: if the $d$-tuple consists of
points $p_1,\dots,p_s$, with multiplicities $m_i$ (so that
$\sum_{i=1}^s m_i=d$), and one puts $m^{(2)}=\sum m_i^2$,
$m^{(3)}=\sum m_i^3$, then the predegree of its orbit-closure equals
$$d^3-3d m^{(2)} +2m^{(3)},$$
so it depends only on $d$ and on $m^{(2)}$, $m^{(3)}$ (one should
compare this result with the degree of the orbit-closure of a smooth
plane curve as computed in Theorem III(B)). Details of this
computation, together with a discussion of the boundary and of
1-dimensional families of orbits, and multiplicity results, can be
found in \cite{Aluffi-Faber}. We should point out that in this case
the degree can also be computed by using simple combinatorics.

Finally, to attract the attention of people working in representation
theory, we remark that we deal here with the orbits of {\sl general\/}
vectors in {\sl one\/} of the standard representations of {\sl one\/}
of the classical algebraic groups. Can these questions be approached
in a more general context? From this point of view, a whole lot of work
remains to be done.

{\smc Acknowledgements.} We thank the University of Chicago, Oklahoma
State University, the Universit\"at Erlangen-N\"urnberg and the
Max-Planck-Institut f\"ur Mathematik for hospitality at various stages
of this work. Discussions with many of our colleagues at these and
other institutions were helpful. We also acknowledge partial support
by the N.S.F. Finally, we are grateful to F. Schubert for writing his
wonderful Fantasia in F minor for piano duet.

\heading \S 2. A blow-up construction\endheading
In this section we construct a smooth projective variety surjecting
onto the orbit closure $\ocbar$ of a {\sl smooth\/} plane curve
$C\in\P^N=\P^\frac{d(d+3)}2$, where $d\ge 3$. As we will see, the
construction depends essentially on the number and type of {\sl
flexes\/} of $C$.

Fix coordinates $(x_0:x_1:x_2)$ of $\P^2$, and assume the degree-$d$
curve $C$ has equation
$$F(x_0,x_1,x_2)=0\quad.$$
Consider the projective space $\P^8=\P\text{Hom}(\C^3,\C^3)$ of
(homogeneous) $3\times 3$ matrices $\alpha=(\alpha_{ij})_{i,j=0,1,2}$.
So $\P^8$ is a compactification of
$\PGL(3)=\{\alpha\in\P^8:\det\alpha\ne 0\}$. To ease notations, in
this section we will refer to a point in $\P^8$ and to any $3\times 3$
matrix representing it by the same term; in the same vein, for
$\alpha\in\P^8$ we will call `$\ker\alpha$' the linear subspace of
$\P^2$ on which the map determined by $\alpha$ is not defined,
`$\im\alpha$' will be the image of this map, and the rank
`$\rk\alpha$' of $\alpha$ will be $1+\dim(\im\alpha)$. So:
$$\alpha\in\PGL(3)\iff\ker\alpha=\emptyset\iff\im\alpha=\P^2
\iff\rk\alpha=3.$$

The curve $C$ determines a rational map
$$c: \P^8 \dashrightarrow \P^N$$
as follows: for $\alpha\in\P^8$, let $c(\alpha)$ be the curve defined
by the degree-$d$ polynomial equation $F(\alpha(x_0,x_1,x_2))=0$. So
$c(\alpha)$ is defined as long as $F(\alpha(x_0,x_1,x_2))$ doesn't
vanish identically; i\.e\., precisely if $\im\alpha\not\subset C$.

If $\alpha\in\PGL(3)$, then $c(\alpha)$ is the translate of $C$ by
$\alpha$; therefore, $c(\PGL(3))$ is just the orbit $\oc$ of $C$ in
$\P^N$ for the natural action of $\PGL(3)$.

As an alternative description for the map $c$, consider for any point
$p\in\P^2$ the equation
$$F(\alpha(p))=0\quad.$$
As an equation `in $p$', this defines the translate $c(\alpha)$; as an
equation `in $\alpha$' this defines the hypersurface of $\P^8$
consisting of all $\alpha$ that map $p$ to a point of $C$. We will
call these hypersurfaces, that will play an important role in our
discussion, `point-conditions'. The rational map defined above is
clearly the map defined by the linear system generated by the
point-conditions on $\P^8$.

Our task here is to resolve the indeterminacies of the map $c: \P^8
\dashrightarrow \P^N$, by a sequence of blow-ups at smooth centers: we
will get a smooth projective variety $\Til V$ filling a commutative
diagram
$$\CD
\PGL(3) \CDsubset \Til V @>\Til c>> \P^N \\
@| @V\pi VV @| \\
\PGL(3) \CDsubset \P^8 \CDdashright c{} \P^N
\endCD$$
The image of $\Til V$ in $\P^N$ by $\Til c$ will then be the orbit
closure $\ocbar$. In \S 3 we will use $\Til c$ to pull-back questions
about $\ocbar$ to $\Til V$; the explicit description of $\Til V$
obtained in this section will enable us to answer these questions.

The plan is to blow-up the support of the base locus of $c$; we will
get a variety $V_1$ and a rational map $c_1:V_1 \dashrightarrow \P^N$.
We will then blow-up the support of the base locus of $c_1$, getting a
variety $V_2$ and a rational map $c_2:V_2 \dashrightarrow \P^N$; in
the case we are considering here (i\.e\., the curve $C$ is smooth to
start with), repeating this process yields eventually a variety $\Til
V$ as above. The support of the first base locus is in fact a copy of
$\P^2 \times C$ in $\P^8$ (see \S2.1); if $(k,q)\in\P^2\times C$, and
$c_i$ denotes the map obtained at the $i$-th stage, we will find that
{\sl $c_i$ still has indeterminacies over $(k,q)$ if and only if the
tangent line to $C$ at $q$ intersects $C$ at $q$ with multiplicity $>
i$.\/} So, for example, if $C$ has only simple flexes then the map
$c_3$ is regular (Proposition 2.9); and in general the number of
blow-ups needed equals the highest possible multiplicity of
intersection of a line with $C$.

We should point out that (even for smooth $C$) this is not the only
way to construct a variety $\Til V$ as above: in fact, a different
sequence of blow-ups is the one that seems to generalize naturally to
approach the same problem for singular $C$.\vskip 12pt

\subheading{\S 2.1. The first blow-up} The set of rank-1 matrices in
$\P^8$ is the image of the Segre embedding
$$\cP^2\times\P^2 \hookrightarrow \P^8$$
given in coordinates by
$$\left( (k_0:k_1:k_2),(q_0:q_1:q_2) \right) \mapsto
\left( \matrix
k_0q_0 & k_1q_0 & k_2q_0 \\
k_0q_1 & k_1q_1 & k_2q_1 \\
k_0q_2 & k_1q_2 & k_2q_2
\endmatrix\right)$$
where $k_0x_0+k_1x_1+k_2x_2=0$ is the kernel of the matrix, and
$(q_0:q_1:q_2)$ is its image. Intrinsically, this is just the map
induced from the map
$$\align
\check{\C}^3 \oplus \C^3 &\rightarrow \check{\C}^3\otimes\C^3=
\text{Hom}(\C^3,\C^3)\\
(f,u)&\mapsto f\otimes u
\endalign$$
We have already observed that the map $c:\P^8 \dashrightarrow \P^N$ is
{\sl not\/} defined at $\alpha\in \P^8$ precisely when
$\im\alpha\subset C$; if $C$ is smooth (therefore irreducible), this
means that the image of $\alpha$ is a point of $C$. Therefore:

\noindent\sl the support of the base locus of $c$ is the image of
$\,\cP^2\times C$ in $\,\P^8$ via the Segre embedding identifying
$\cP^2\times\P^2$ with the set of rank-1 matrices.\rm

In particular, the support of the base locus of $c$ is {\sl smooth\/},
since $C$ is. We let then $B=\cP^2\times C$, and we let $V_1 @>\pi_1
>> \P^8$ be the blow-up of $\P^8$ along $B$. Since
$B\cap\PGL(3)=\emptyset$, $V_1$ contains a dense open set which we can
identify with $\PGL(3)$. Also, the linear system generated by the
proper transforms in $V_1$ of the point-conditions (which we will call
`point-conditions in $V_1$'), defines a rational map $c_1:V_1
\dashrightarrow \P^N$ making the diagram
$$\CD
\PGL(3) \CDsubset V_1 \CDdashright {c_1}{} \P^N \\
@| @V\pi_1 VV @| \\
\PGL(3) \CDsubset \P^8 \CDdashright c{} \P^N
\endCD$$
commutative. The exceptional divisor $E_1$ in $V_1$ is the
projectivized normal bundle of $B$ in $\P^8$: $E_1=\P(N_B\P^8)$. We
will show now that {\sl the base locus of $c_1$ is supported on a
$\P^1$-subbundle of $E_1$ over $B$.\/}

Let $(k,q)$ be a point of $B=\cP^2\times C$: i\.e\., a rank-1
$\alpha\in\P^8$ with $\ker\alpha=k$, $\im\alpha=q\in C$. Also, let
$\ell$ be the line tangent to $C$ at $q$, let $p$ be a point of
$\P^2$, and denote by $P$ the point-condition in $\P^8$ corresponding
to $p$.
\proclaim{Lemma 2.1} (i) The tangent space to $B$ at $(k,q)$ consists
of all $\phy\in\P^8$ such that $\im\phy\subset \ell$ and $\phy(k)
\subset q$.

(ii) $P$ is non-singular at $(k,q)$, and the tangent space to $P$ at
$(k,q)$ consists of all $\phy\in\P^8$ such that $\phy(p)\subset\ell$.
\endproclaim
We are using our notations rather freely here. For example, in (i)
$\alpha =(k,q)$ is in the tangent space since $\alpha(k)=\emptyset$
(as $\alpha$ is not defined along $k$).
\demo{Proof} (i) The tangent space to $B$ at $(k,q)$ is spanned by the
plane $\{(k',q)\in B:k'\in\cP^2\}=\{\phy\in\P^8:\im\phy=q\}$ and by
the line $\{(k,q')\in
B:q'\in\ell\}=\{\phy\in\P^8:\ker\phy=k,\im\phy\in\ell\}$. Both these
subspaces of $\P^8$ are contained in
$\{\phy\in\P^8:\im\phy\subset\ell, \phy(k)\subset q\}$; since this
latter has clearly dimension 3, we are done.

(ii) For $\alpha=(k,q)$ and $\phy\in\P^8$ consider the line
$\alpha+\phy\,t$.  Restricting the equation for $P$ to this line gives
the polynomial equation in $t$
$$\gather
F((\alpha+\phy\,t)(p))=0\quad,\text{ i\.e\.}\\
F(\alpha(p))+\sum_i\left(\frac{\partial F}{\partial x_i}\right)_{\alpha(p)}
\phy_i(p)\,t+\dots=0
\endgather$$
(where $\phy_i(p)$ denotes the $i$-th coordinate of $\phy(p)$).

$F(\alpha(p))=0$ since $\im\alpha=q\in C$; the line is tangent to $P$
at $\alpha$ when the linear term also vanishes, i\.e\. if $\sum_i
(\partial F/ \partial x_i)_q \phy_i(p)=0$. This says precisely
$\phy(p)\subset\ell$, as claimed.

$P$ is non-singular at $\alpha$ because any $\phy$ {\sl not\/}
satisfying the condition $\phy(p)\subset\ell$ gives a line
$\alpha+\phy\,t$ intersecting $P$ with multiplicity 1 at $\alpha$, by
the above computation.\qed\enddemo

With the same notations, the tangent space to $\cP^2\times\P^2$ at
$\alpha$ consists of all $\phy$ with
$\phy(\ker\alpha)\subset\im\alpha$ (intrinsically, all transformations
$\phy$ inducing a map coim$\,\alpha @>>>\,$coker$\,\alpha$).

The set of all $\phy$ such that $\im\phy\subset \ell$ forms (for any
$\alpha$) a 5-dimensional space containing the tangent space to $B$ at
$\alpha$, and therefore determines a 2-dimensional subspace of the
fiber of $N_B\P^8$ over $\alpha$. As $\alpha$ moves in $B$ we get a
rank-2 subbundle of $N_B\P^8$, and hence a $\P^1$-subbundle of
$E_1=\P(N_B\P^8)$, which we denote $B_1$. Notice that $B_1$ is
non-singular, as a $\P^1$-bundle over the non-singular $B$.
\proclaim{Proposition 2.2} The base locus of the map $c_1:V_1
\dashrightarrow \P^N$ is supported on $B_1$.\endproclaim
\demo{Proof} Since $c_1$ is defined by the linear system generated by
all point-conditions in $V_1$, we simply need to show that the
intersection of all point-conditions in $V_1$ is set-theoretically
$B_1$. This assertion can be checked fiberwise over $\alpha=(k,q)\in
B$; so all we need to observe is that the intersection of the tangent
spaces to all point-conditions at $\alpha$ consists (by Lemma 2.1
(ii)) of the $\phy\in\P^8$ such that $\phy(p)\subset\ell$ for {\sl
all\/} $p$; i\.e\., the 5-dimensional space used above to define
$B_1$.\qed\enddemo

If $P_1^{(p)}$ denotes the point-condition in $V_1$ corresponding to
$p\in\P^2$, we have just shown $\bigcap_{p\in\P^2} P_1^{(p)}$ is
supported on $B_1$. The proof says a little more:
\proclaim{Remark 2.3} $\bigcap_{p\in\P^2} P_1^{(p)}\cap E_1=B_1$
(scheme-theoretically).\endproclaim
Indeed on each fiber of $E_1$ (say over $\alpha\in B$) the fiber of
$B_1$, a linear subspace, is cut out by the fibers of the
$P_1^{(p)}\cap E_1$, linear subspaces themselves; and the situation
clearly globalizes as $\alpha$ moves in $B$.\vskip 12pt

\subheading{\S 2.2. The second blow-up} Let $V_2 @>\pi_2 >> V_1$ be
the blow-up of $V_1$ along $B_1$. The new exceptional divisor is
$E_2=\P(N_{B_1} V_1)$; call `point-conditions in $V_2$' the proper
transforms of the point-conditions of $V_1$. The linear system
generated by the point-conditions defines a rational map $c_2: V_2
\dashrightarrow \P^N$; again, we obtain a diagram
$$\CD
\PGL(3) \CDsubset V_2 \CDdashright {c_2}{} \P^N \\
@| @V\pi_2\circ\pi_1 VV @| \\
\PGL(3) \CDsubset \P^8 \CDdashright c{} \P^N
\endCD$$
and we proceed to determine the support of the base locus of $c_2$.

Let $\Til E_1$ be the proper transform of $E_1$ in $V_2$. Then
\proclaim{Lemma 2.4} The base locus of $c_2$ is disjoint from $\Til E_1$.
\endproclaim
\demo{Proof} This is basically a reformulation of Remark 2.3: $\Til E_1$ is
the blow-up of $E_1$ along $B_1$, and $B_1$ is cut out scheme-theoretically by
the intersections of $E_1$ with the point-conditions of $V_1$. So the
intersection of the point-conditions in $V_2$ must be empty along $\Til E_1$,
which is the claim.\qed\enddemo

Lemma 2.4 reduces the determination of the support of the base locus
of $c_2$ to a computation in $\P^8$. Denote by $\Cal B$ the
scheme-theoretic intersection of the point-conditions in $\P^8$, so
the support of $\Cal B$ is $B$. For $\alpha\in B$, let $th_\alpha(\Cal
B)$ be the maximum length of the intersection with $\Cal B$ of the
germ of a smooth curve centered at $\alpha$ and transversal to $B$
(the `thickness' of $\Cal B$ at $\alpha$, in the terminology of
\cite{Aluffi2}).
\proclaim{Lemma 2.5} The base locus of $c_2$ is disjoint from $(\pi_2\circ
\pi_1)^{-1}\alpha$ if $th_\alpha(\Cal B)\le 2$.\endproclaim
\demo{Proof} The base locus of $c_2$ is the intersection of all
point-conditions in $V_2$, i\.e\. the set of all directions normal to $B_1$
and tangent to all point-conditions in $V_1$. Let then $\gamma(t)$ be a smooth
curve germ centered at a point of $B_1$ above $\alpha$, transversal to $B_1$,
and tangent to all point-conditions in $V_1$. By Lemma 2.4, $\gamma$ is
transversal to $E_1$; therefore $\pi_1(\gamma(t))$ is a {\sl smooth\/} curve
germ centered at $\alpha$ and {\sl transversal\/} to $B$. Since $\gamma(t)$
intersects all point-conditions in $V_1$ with multiplicity 2 or more,
$\pi_1(\gamma(t))$ must intersect all point-conditions in $\P^8$ with
multiplicity 3 or more; $\Cal B$ is the intersection of all point-conditions
in $\P^8$, so this forces $th_\alpha(\Cal B)\ge 3$.\qed\enddemo

Now the key computation is
\proclaim{Lemma 2.6} If $\alpha=(k,q)\in B$, and $\ell$ is the line tangent to
$C$ at $q$, then $th_\alpha(\Cal B)$ equals the intersection multiplicity of
$\ell$ and $C$ at $q$.\endproclaim
\demo{Proof} Let $m$ be the intersection multiplicity of $\ell$ and $C$ at
$q$. To show $th_\alpha(\Cal B)\ge m$, we just have to produce a curve normal
to $B$ and intersecting all point-conditions with multiplicity at least $m$ at
$\alpha$; such is the line $\alpha+\phy\,t$, with $\phy\in\P^8$ such that
$\im\phy=\ell$ and $\phy(k)\ne q$. Indeed, the last condition guarantees
normality (Lemma 2.1 (i)); and, for general $p$, $q=\alpha(p)$ and $\phy(p)$
span $\ell$: so $F((\alpha+\phy\,t)(p))$ is just the restriction of $F$ to a
parametrization of $\ell$, and it must vanish exactly $m$ times at $t=0$.
Notice that these directions are precisely those defining $B_1$.

To show $th_\alpha(\Cal B)\le m$, let $\gamma(t)$ be any smooth curve germ
normal to $B$ and centered at $\alpha$; we have to show that $\gamma$
intersects some point-condition with multiplicity $\le m$ at $\alpha$. In an
affine open of $\P^8$ containing $\alpha$, write
$$\gamma(t)=\alpha+\phy\,t+\dots\quad.$$
The equation for the point-condition corresponding to $p$ restricts on
$\gamma$ to
$$F((\alpha+\phy\,t+\dots)(p))= F(\alpha(p))+\sum_i\left(\frac{\partial F}
{\partial x_i}\right)_{\alpha(p)} \phy_i(p)\,t+\dots=0\quad,$$
where $\phy_i(p)$ denotes the $i$-th coordinate of $\phy(p)$. The coefficient
of $t^m$ in this expansion is
$$\multline
\frac 1{m!}\sum_{i_1,\dots,i_m}\left(\frac{\partial^mF}{\partial x_{i_1}
\cdots\partial x_{i_m}}\right)_{\alpha(p)} \phy_{i_1}(p)\cdots \phy_{i_m}(p)\\
+\text{ terms involving derivatives of lower order,}
\endmultline\tag*$$
and to conclude the proof we have to show that for some $p$ this term
doesn't vanish.

To see this, observe that since $\ell$ and $C$ intersect with multiplicity
exactly $m$ at $q$, then the form
$$\sum_{i_1,\dots,i_m}\left(\frac{\partial^mF}{\partial x_{i_1} \cdots\partial
x_{i_m}} \right)_{\alpha(p)}x_{i_1}\cdots x_{i_m}$$
doesn't vanish identically on $\ell$; since $\phy(\ker\alpha)\not\subset q$
($\gamma$ is normal to $B$), this implies that the summand
$$\frac 1{m!}\sum_{i_1,\dots,i_m}\left(\frac{\partial^mF}{\partial x_{i_1}
\cdots\partial x_{i_m}}\right)_{\alpha(p)} \phy_{i_1}(p)\cdots \phy_{i_m}(p)$$
vanishes exactly $d-m$ times along the line $k=\ker\alpha$ (as a
function of $p$). But since all the other summands in (*) involve
derivatives of order $<m$, they vanish with order $>d-m$ along $k$.
Therefore the order of vanishing of (*) along $k$ must be exactly
$d-m$, and in particular (*) can't be identically 0, as we
claimed.\qed\enddemo

We adopt the following convention:
\proclaim{Definition} A point $q$ of $C$ is a `flex of order $r$' if
the line tangent to $C$ at $q$ intersects $C$ at $q$ with multiplicity
$r+2$. We will say that $q$ is a `flex' of $C$ if $r\ge 1$, and that
$q$ is a `simple flex' if $r=1$.\endproclaim

Now we observe that there is a section $s:B_1 @>>> E_2$: for $\alpha_1\in
B_1$, let $\alpha=\pi_1(\alpha_1)\in B$, say $\alpha=(k,q)$, and let $\ell$ be
the line tangent to $C$ at $q$. By the construction of $B_1$, there is a
matrix $\phy\in\P^8$ with $\im\phy\subset\ell$ such that $\alpha_1$ is the
intersection of $E_1$ and the proper transform of the line $\alpha+\phy\,t$ in
$V_1$; now let $s(\alpha_1)$ be the intersection of $E_2$ and the proper
transform of the line $\alpha+\phy\,t$ in $V_2$ (it is clear that
$s(\alpha_1)$ does not depend on the specific $\phy$ chosen to represent
$\alpha_1$).

Let $B_2$ be the image via $s$ of the set $\{\alpha_1\in B_1:q \text{
is a flex of $C$}\}$. Thus $B_2$ consists of a number of smooth
three-dimensional components, one for each flex of $C$: each
component maps isomorphically to a $\P^1$-bundle over one of the
planes $\{(k,q)\in B:q \text{ is a flex of $C$}\}$.

\proclaim{Proposition 2.7} The base locus of the map $c_2: V_2 \dashrightarrow
\P^N$ is supported on $B_2$.\endproclaim

\demo{Proof} Let $\alpha_1\in B_1$, and $\alpha=(k,q)$ the image of
$\alpha_1$ in $B$, as above.  Consider the intersection of the base
locus of $c_2$ with the fiber $\pi_2^{-1}(\alpha_1)\cong \P^3$. By
Lemma 2.5 and 2.6 this is {\sl empty\/} if $q$ is not a flex of $C$;
even if $q$ is a flex of $C$, the intersection is a linear subspace of
$\P^3$ missing a $\P^2$ (by Lemma 2.4), thus it consists of at most
{\sl one point.\/} Thus all we have to show is that $s(\alpha_1)$ {\sl
is\/} in the base locus of $c_2$ if $q$ is a flex of $C$ (of order
$r\ge 1$). But, as observed in the proof of Lemma 2.6, the line
$\alpha+\phy\,t$ determining $\alpha_1$ intersects each
point-condition in $\P^8$ with multiplicity at least $r+2\ge 3$;
therefore the proper transform of $\alpha+\phy\, t$ is {\sl tangent\/}
to all point-conditions in $V_1$, and it follows that $s(\alpha_1)\in$
all point-conditions in $V_2$, as needed.\qed\enddemo
\vskip 12pt

\subheading{\S 2.3. The third blow-up} Let $V_3 @>\pi_3>> V_2$ be the blow-up
of $V_2$ along $B_2$. The new exceptional divisor is $E_3$; the
`point-conditions of $V_3$' are the proper transforms of the point-conditions
of $V_2$. The linear system generated by the point-conditions defines a
rational map $c_3: V_3 \dashrightarrow \P^N$, making the diagram
$$\CD
\PGL(3) \CDsubset V_3 \CDdashright {c_3}{} \P^N \\
@| @V\pi_3\circ\pi_2\circ\pi_1 VV @| \\
\PGL(3) \CDsubset \P^8 \CDdashright c{} \P^N
\endCD$$
commute. We will show now that $c_3$ is a {\sl regular map\/} if all the
flexes of $C$ are simple, so that in this case $V_3$ is the variety we are
looking for. For each flex of order $> 1$, we will find a four-dimensional
component in the base locus of $c_3$, and more blow-ups will be needed.

Call $\Cal B_2$ the scheme-theoretic intersection of the point-conditions in
$V_2$, so $\Cal B_2$ is supported on $B_2$. For $\alpha_2\in B_2$, define the
thickness $th_{\alpha_2}(\Cal B_2)$ of $\Cal B_2$ at $\alpha_2$ as we did
above for $th_\alpha(\Cal B)$. Also, let $\alpha=(k,q)$ be the image of
$\alpha_2$ in $B$. With these notations:
\proclaim{Lemma 2.8} If $q$ is an flex of order $r$ of $C$, then
$th_{\alpha_2}(\Cal B_2)=r$.\endproclaim

\demo{Proof} We have to show that if $\gamma(t)$ is a smooth curve germ in
$V_2$, centered at $\alpha_2$ and transversal to $B_2$, then the maximum
length of the intersection of $\Cal B_2$ and $\gamma$ at $t=0$ is precisely
$r$.

Suppose first that $\gamma$ is transversal to $E_2$: then, as argued
in the proof of Lemma 2.5, the image of $\gamma$ in $\P^8$ is a smooth
curve germ centered at $\alpha$ and transversal to $B$: by Lemma 2.6,
the length of the intersection of $\Cal B$ and this curve is at most
$r+2$; it follows that the maximum length of the intersection of $\Cal
B_2$ and such $\gamma$'s is indeed $r$ (attained for example by the
proper transform of $\alpha+\phy\,t$, with $\phy$ as in the proof of
Lemma 2.6).

Thus we may assume that $\gamma$ is tangent to $E_2$, and we have to
show that
\proclaim{Claim} $\Cal B_2\cap \gamma(t)$ vanishes at most $r$ times at $t=0$.
\endproclaim
This is a lengthy but straightforward coordinate computation, which we
leave to the reader.
The outcome is that the maximum length is $r$, and it
is attained in the direction normal to $B_2$ in the section $s(B_1)\subset
E_2$ defined in \S 2.2.
\qed\enddemo

The next results are now easy consequences.
\proclaim{Proposition 2.9} If all flexes of $C$ are simple, then the map
$c_3:V_3 \dashrightarrow \P^N$ is regular.\endproclaim

\demo{Proof} We have to show that $c_3$ has no base locus, i\.e\. that the
intersection of all point-conditions in $V_3$ is {\sl empty\/}. But a point in
the intersection of all point-conditions in $V_3$ would determine a direction
normal to $B_2$ and {\sl tangent\/} to all point-conditions in $V_2$; the
thickness of $\Cal B_2$ would then be $\ge 2$ at some point. By Lemma 2.8, if
all flexes of $C$ are simple (i\.e\., of order 1) the thickness of $\Cal B_2$
is precisely 1 everywhere on $B_2$, so this can't happen.\qed\enddemo

By Proposition 2.9, we are done in the case when $C$ has only simple flexes:
$V_3$ is the variety $\Til V$ we meant to construct. We will show now that for
each flex of $C$ of order $r>1$, the base locus of $c_3$ has a smooth
four-dimensional connected component.

Let $\alpha_2\in B_2$, mapping to $\alpha=(k,q)$ in $B$, and assume
$q$ is a flex of $C$ of order $r>1$. $B_2$ is 3-dimensional, so the
fiber $\pi^{-1}_3(\alpha_2)$ of $E_3=\P(N_{B_2}V_2)$ over $\alpha_2$
is a $\P^4$.  We have two special points in this $\P^4$, namely the
point determined by the proper transform of the line $\alpha+\phy\,t$
used in \S 2.2 to define $s$, and the direction normal to $B_2$ in the
section $s(B_1)$. We have seen in the proof of Lemma 2.8 that the
length of the intersection of these directions with $\Cal B_2$ is
exactly $r$; also, these points are {\sl distinct\/} for all
$\alpha_2$ (since one of them corresponds to a direction contained in
$E_2$, while the other comes from a direction transversal to $E_2$),
so they determine a $\P^1$ in the fiber $\pi^{-1}_3(\alpha_2)$. As
$\alpha_2$ moves in the component of $B_2$ over $q$, this $\P^1$
traces a $\P^1$-bundle over that component, a smooth four-dimensional
subvariety $\bq_3$ of $E_3$. Call $B_3$ the union of all these
(disjoint) subvarieties of $E_3$, arising from non-simple flexes of
$C$.

\proclaim{Proposition 2.10} The base locus of the map $c_3: V_3
\dashrightarrow \P^N$ is supported on $B_3$.\endproclaim
\demo{Proof} The argument here is somewhat analogous to the argument
in the proof of 2.7. We have to show that in each fiber
$\pi_3^{-1}(\alpha_2)\cong\P^4$ as above, the intersection of all
point-conditions is supported on the specified $\P^1$. Observe that each
point-condition determines a hyperplane in this $\P^4$, so that the
intersection of the base locus of $c_3$ with $\pi_3^{-1}(\alpha_2)$ must be a
linear subspace of this $\P^4$. Secondly, for the same reason, no directions
tangent to the fiber of $E_2$ containing $\alpha_2$ can be tangent to all
point-conditions in $V_2$.
The fibers of $E_2$ are three-dimensional and
transversal to $B_2$, thus this shows that the base locus of $c_3$ must miss a
$\P^2$ in the fiber $\pi_3^{-1}(\alpha_2)$. Thus, the intersection of the base
locus of $c_3$ with $\pi_3^{-1}(\alpha_2)$ can consist of at most a $\P^1$.

Therefore, we just have to show that the two points of $\pi_3^{-1}(\alpha_2)$
used in the construction of $B_3$ are contained in all point-conditions of
$V_3$; or, equivalently, the two directions in $V_2$ used to define these
points are {\sl tangent\/} to all point-conditions in $V_2$. But this is
precisely the result of the computation in the proof of Lemma 2.8: the length
of the intersection of these curves with all point-conditions is $r\ge 2$.
\qed\enddemo\vskip 12pt

\subheading{\S 2.4 Further blow-ups} As we have seen in \S 2.3, each
non-simple flex $q$ of $C$ gives rise to a smooth four-dimensional
component of the support $B_3$ of the base locus of $c_3$; and $B_3$
is the union of all such components. The plan is still to blow-up the
support of the base-locus; since the components are disjoint, we can
concentrate on a specific one: say $\bq_3$, corresponding to a flex
$q$ of $C$ of order $r\ge 2$.

Let $\vq_3$ be the complement of all components of $B_3$ {\sl other than\/}
$\bq_3$ in $V_3$. Let $\vq_4 @>>> \vq_3$ be the blow-up of $\vq_3$ along
$\bq_3$; again, the proper transforms in $\vq_4$ of the point-conditions
define a map $\cq_4: \vq_4 \dashrightarrow \P^N$. The base locus of
$\cq_4$ might have components over $\bq_3$, whose union we denote
$\bq_4$; in this case, we will let $\vq_5$ be the blow-up of $\vq_4$
along $\bq_4$. Iterating this process we get a tower of varieties and maps:
$$\CD
\vdots       @.      \vdots                 @.          \vdots \\
\bq_{i+1} \CDsubset \vq_{i+1} \CDdashright{\cq_{i+1}}{} \P^N \\
@VVV                @VVV                                  @|   \\
\bq_i     \CDsubset \vq_i     \CDdashright{\cq_i}{}     \P^N \\
@VVV                @VVV                                  @|   \\
\vdots       @.      \vdots                 @.          \vdots \\
\bq_3     \CDsubset \vq_3     \CDdashright{c_3}{}       \P^N
\endCD$$
where, inductively for $i\ge 4$: $\vq_i @>>> \vq_{i-1}$ is the blow-up
of $\vq_{i-1}$ along $\bq_{i-1}$; $\cq_i:\vq_i \dashrightarrow \P^N$
is defined by the proper transforms in $\vq_i$ of the point-conditions
(i\.e\., the `point-conditions in $\vq_i$');
and (for $i\ge 3$) $\bq_i$ is the support of
the intersection $\Bq_i$ of the point-conditions in $\vq_i$ (i\.e\., the base
locus of $\cq_i$). Also, for $i\ge 3$ let $\Eq_i$ be the exceptional divisor
in $\vq_i$, and let $\Til E^{(q)}_i$ be the proper transform of $\Eq_i$ in
$V^{(q)}_{i+1}$.

\proclaim{Lemma 2.11} If $q$ is a flex of order $r\ge 2$, then
for $3\le i\le r+1$:

$(1)_i$: $\vq_i$ is non-singular

$(2)_i$: the composition map $\bq_i @>>> \bq_3$ is an isomorphism

$(3)_i$: the thickness of $\Bq_i$ is $r+2-i$ at each point of $\bq_i$

$(4)_i$: $\bq_{i+1}\cap \Til E^{(q)}_i=\emptyset$\endproclaim

\demo{Proof} We have $(1)_3$, $(2)_3$ trivially, and $(3)_2$ by Lemma 2.8.
Also, since $B_3$ is cut out by linear spaces in each fiber of $E_3$,
we have $(4)_3$. Now we will show that:
\proclaim{Claim} For $4\le i\le r+1$, $(1)_{i-1}$, $(2)_{i-1}$, $(3)_{i-2}$
and $(4)_{i-1}$ imply $(1)_i$, $(2)_i$, $(3)_{i-1}$, and $(4)_i$.\endproclaim
Also, we will show that $(3)_r$, $(4)_{r+1}$ imply $(3)_{r+1}$: this will
prove the statement.
\demo{Proof of the Claim} In this proof we will drop the $\,^{(q)}$ notation,
to ease the exposition. $V_i$ is then the blow-up of $V_{i-1}$ along
$B_{i-1}$, and these are both non-singular by $(1)_{i-1}$, $(2)_{i-1}$: so
$V_i$ must also be non-singular, giving $(1)_i$.

Next, compute the thickness of $\Cal B_{i-1}$: let $\gamma(t)$ be any smooth
curve germ transversal to $B_{i-1}$ and centered at any $\alpha_{i-1}\in
B_{i-1}$. If $\gamma$ is tangent to $E_{i-1}$, then by $(4)_{i-1}$ its proper
transform will miss the general point-condition in $V_i$: i\.e\., the length
of the intersection of $\gamma(t)$ with $\Cal B_{i-1}$ at $t=0$ is 1. If
$\gamma$ is transversal to $E_{i-1}$ (and $B_{i-1}$), then $\gamma$ maps down
to a {\sl smooth\/} curve germ $\gamma_*$ centered at a point of $B_{i-2}$ and
transversal to $B_{i-2}$. By $(3)_{i-2}$, the intersection of $\gamma_*$ with
the point-conditions in $V_{i-2}$ has length at most $r-i+4$: it follows that
the intersection of $\gamma$ with the point-conditions in $V_{i-1}$ has length
at most $r-i+3\ge 2$ (since $i\le r+1$). Therefore the thickness of
$\Cal B_{i-1}$ at $\alpha_{i-1}$ is $r-i+3$, which gives $(3)_{i-1}$.

For $(2)_i$, look at the intersection of $B_i$ with the fiber of $E_i$ over
an arbitrary $\alpha_{i-1}\in B_{i-1}$. First we argue this can't be empty:
indeed, $th_{\alpha_{i-1}}(\Cal B_{i-1})=r-i+3\ge 2$, so through every
$\alpha_{i-1}$ in $B_{i-1}$ there are directions tangent to all
point-conditions in $V_{i-1}$. To get $(2)_i$, we need to show that the fiber
of $B_i$ over $\alpha_{i-1}$ consists (scheme-theoretically) of a simple
point. But this is the intersection of $B_i$ with the fiber of $E_i$
($\cong\P^3$) over $\alpha_{i-1}$, thus a nonempty intersection of linear
subspaces in $\P^3$ missing a hyperplane (by $(4)_{i-1}$): precisely a point,
as needed for $(2)_i$.

Finally, we need $(4)_i$. Once more observe that $\Cal B_i$ intersects each
fiber of $E_i$ in an intersection of linear spaces: thus there are no
directions in the fibers of $E_i$ and tangent to all point-conditions in
$V_i$. This says that $B_{i+1}$ must avoid the proper transforms in $V_{i+1}$
of all fibers of $E_i$, and therefore $\Til E_i$, giving $(4)_i$.\enddemo

This proves the Claim. The only case not covered yet is $(3)_{r+1}$: to obtain
this and conclude the proof of 2.11, apply the same argument as above to
$(3)_r$, $(4)_{r+1}$.\qed\enddemo
Lemma 2.11 describes the sequence of blow-ups over $V_3$ that takes care
of a specific flex $q$ on $C$ of order $r\ge 2$. The case $i=r+1$ of the
statement says that the variety $\vq_{r+1}$ is non-singular, and the base
locus of the map $\cq_{r+1}: \vq_{r+1}\dashrightarrow \P^N$ is supported on a
variety $\bq_{r+1}$ isomorphic to $\bq_3$; moreover, for all $\alpha_{r+1}\in
\bq_{r+1}$, we got $th_{\alpha_{r+1}}(\Cal B_{r+1})=1$. Let then $\vq_{r+2}
@>>> \vq_{r+1}$ be the blow-up of $\vq_{r+1}$ along $\bq_{r+1}$, and denote by
$\cq_{r+2}$ the rational map $\vq_{r+2} \dashrightarrow \P^N$ defined by the
point-conditions in $\vq_{r+2}$. Then $\vq_{r+2}$ is clearly non-singular, and
\proclaim{Corollary 2.12} $\cq_{r+2}$ is a regular map.\endproclaim
\demo{Proof} Indeed, the point-conditions in $\vq_{r+2}$ cannot intersect
anywhere along $\Eq_{r+2}$: if they did, any intersection point would
correspond to a direction normal to $\bq_{r+1}$ and tangent to all
point-conditions in $\vq_{r+1}$, and the thickness of $\Bq_{r+1}$ would be
$\ge 2$, in contradiction with Lemma 2.11.\qed\enddemo

By this last result, the sequence of $r-1$ blow-ups over $V_3$ just described
resolves the indeterminacies of $c_3: V_3 \dashrightarrow \P^N$ over the
component $\bq_3$ of $B_3$. To resolve {\sl all\/} indeterminacies of $c_3$,
we just have to apply the construction simultaneously to all components of
$B_3$: build the sequence
$$\CD
\vdots       @.      \vdots                 @.          \vdots \\
B_{i+1} \CDsubset V_{i+1} \CDdashright{c_{i+1}}{} \P^N \\
@VVV                @VVV                                  @|   \\
B_i     \CDsubset V_i     \CDdashright{c_i}{}     \P^N \\
@VVV                @VVV                                  @|   \\
\vdots       @.      \vdots                 @.          \vdots \\
B_3     \CDsubset V_3     \CDdashright{c_3}{}       \P^N
\endCD$$
where, for $i\ge 4$, $V_i @>>> V_{i-1}$ is the blow-up of $V_{i-1}$ along
$B_{i-1}$, $c_i: V_i \dashrightarrow \P^N$ is defined by the proper transforms
in $V_i$ of the point-conditions,
and $B_i$ is the support of the base locus
of $c_i$. By Lemma 2.11 and Corollary 2.12 all $V_i$'s are non-singular,
and, for each flex $q$ of $C$ of order $r$, $B_i$ has either exactly
one component mapping isomorphically to $\bq_3$ if $i\le r+1$, or {\sl
no\/} component over $\bq_3$ if $i\ge r+2$.

In particular, this construction will stop! If $r$ is the maximum among the
order of the flexes of $C$, let $\Til V=V_{r+2}$, $\Til c=c_{r+2}$, and let
$\pi$ be the composition of the $r+2$ blow-up maps; then we have shown
\proclaim{Theorem II} $\Til c: \Til V @>>> \P^N$ is a regular map, and the
diagram
$$\CD
\PGL(3) \CDsubset \Til V @>\Til c>> \P^N \\
@| @V\pi VV @| \\
\PGL(3) \CDsubset \P^8 \CDdashright c{} \P^N
\endCD$$
commutes.\endproclaim
\noindent which was our objective.

\heading \S 3. The degree of the orbit closure\endheading In this
section we employ the blow-up construction of $\S 2$ to compute the
degree of the orbit closure $\ocbar$ of a smooth plane curve
$C\in\P^N=\P^{\frac{d(d+3)}2}$ with at most finitely many
automorphisms (if $d=3$, we should specify `induced from $\PGL(3)$'.
This will be understood in the following). The degree will depend on
just {\sl six\/} natural numbers: the order of the group of
automorphisms of $C$, the degree $d$ of $C$, and four numbers
encoding information about the number and order of the flexes of~$C$.
In fact, the blow-up construction of \S 2 yields most naturally the
`predegree' of $\oc$:
\proclaim{Definition} The `predegree' of $\oc$ is the 8-fold self-intersection
$\Til P^8$ of the class $\Til P$ of a point-condition in $\Til V$.
\endproclaim
\proclaim{Lemma 3.1} The predegree of $\oc$ equals the product of the degree
of the orbit closure of $C$ by the order $o_C$ of the group of
automorphisms of $C$ induced from $\PGL(3)$.\endproclaim
\demo{Proof} The map $\Til c$ is defined by the linear system generated by the
point-conditions on $\Til V$, so $\Til P$ is the pull-back of the
hyperplane class from $\P^N$. Therefore $\Til P^8$ computes the
pull-back of the intersection of $\Til c(\Til V)=\ocbar$ with 8
hyperplanes of $\P^N$: i\.e\., the product of $\deg(\ocbar)$ by the
degree of the map $\Til c$. This latter equals $o_C$ since, given a
general $c(\alpha)\in\oc$
($\alpha\in\P^8$), the fiber of $c(\alpha)$ consists of all products
$\phy\alpha$, where $\phy$ fixes $C$.
\qed\enddemo

Observe that for the general $C$ of degree $\ge 4$, the predegree of
$\oc$ equals the degree of the orbit closure. Our aim here is to
compute the predegree of $\oc$, by using the construction of $\Til V$
described in \S 2: we will show that this number depends only on $d$
and on the type of the flexes of $C$.\vskip 12pt

Our tool will be a formula relating intersection degrees under blow-ups:
\proclaim{Proposition 3.2} Let $B \overset i \to\hookrightarrow V$ be
non-singular projective varieties, and let $X\subset V$ be a codimension-1
subvariety, smooth along $B$. Let $\Til V$ be the blow-up of $V$ along $B$,
and let $\Til X$ be the proper transform of $X$. Then
$$\int_{\Til V} [\Til X]^{\dim V}=\int_V [X]^{\dim V} - \int_B\frac{([B]+i^*
[X])^{\dim V}}{c(N_BV)}$$\endproclaim
\noindent where $\int_{\Til V}$, etc\. denote the degree of a class in $\Til
V$, etc\., cf\. \cite{Fulton}, Def\. 1.4. Note: we will omit the
$\int$ sign and the class brackets $[\cdot]$ when this doesn't create
ambiguities.
\demo{Proof} This follows from \cite{Aluffi1}, \S2, Theorem II and
Lemma (2), (3).\qed\enddemo

We will compute the predegree of $\oc$ (i\.e\. $\Til P^8$) by applying
Proposition 3.2 to each blow-up in the sequence giving $\Til V$: the
missing ingredients to be obtained at this point are the Chern classes
of the normal bundles of the centers of the blow-ups, and calculations
in their intersection rings.

In the following, $P,P_i,\Til P$ will denote resp\. (the class of)
point-conditions in $V,V_i,\Til V$. The embedding of $B_j$ in $V_j$ is denoted
$i_j$, and $p_{jk}$ will be used for the map $B_j @>>> B_k$ ($p_j$ will be
$p_{jj-1}$ for short). As a general convention, we will omit pull-back
notations unless we fear ambiguity.\vskip 12pt

\subheading{\S 3.1. The first blow-up} The center of the first blow-up is the
variety $B=\cP^2\times C$; the embedding $i: B\hookrightarrow \P^8$ is given
by composition with the Segre embedding:
$$B=\cP^2\times C\subset \cP^2\times\P^2 @>>> \P^8\quad.$$
Call $h,k$ resp\. the hyperplane class in $\P^2,\cP^2$. Our convention on
pull-backs allows us to write $k,h$ for the pull-backs of $k,h$ from the
factors to $\cP^2\times\P^2$, and to $B\subset\cP^2\times\P^2$. Also, since
the Segre embedding is linear on each factor, the hyperplane class of $\P^8$
pulls-back to $k+h$ on $B$.

\proclaim{Lemma 3.3} If $C$ has degree $d$:

(i) In $B$: $k^3=0,k^2h=d,kh^2=0,h^3=0$

(ii) $c(N_B\P^8)=\dfrac{(1+k+h)^9(1+dh)}{(1+k)^3(1+h)^3}$

(iii) $P^8=d^8$; and $P$ pulls-back to $dk+dh$.\endproclaim
\demo{Proof} (i) is immediate.

(ii) $c(N_B\P^8)=c(N_B\cP^2\times\P^2)c(N_{\cP^2\times\P^2}\P^8)$ by the
Whitney formula and the exact sequence of normal bundles. Now, since
$B=\cP^2\times C$, $c(N_B\cP^2\times\P^2)=c(N_C\P^2)=1+dh$. The formula for
$c(N_{\cP^2\times\P^2}\P^8)$ is standard.

(iii) Recall from \S 2 that if $p\in\P^2$, $P$ is the point-condition
corresponding to $p$, and $F(x_0:x_1:x_2)$ is the (degree-$d$) polynomial
defining $C$, then $\alpha\in P\iff F(\alpha(p))=0$: so $P$ is defined by a
degree-$d$ equation in $\P^8$.\qed\enddemo

We have already observed that the point-conditions are non-singular (Lemma 2.1
(ii)), so we are ready for the key computation needed to apply Proposition 3.2
to the first blow-up:
\proclaim{Lemma 3.4}
$$\int_B\frac{(B+i^*P)^8}{c(N_B\P^8)}=d(10d-9)(14d^2-33d+21)$$
\endproclaim
\demo{Proof} By Lemma 3.3, this is
$$\int_{\cP^2\times C}\frac{(1+dk+dh)^8(1+k)^3(1+h)^3}{(1+k+h)^9(1+dh)}
\quad:$$
the statement follows by computing the coefficient of $k^2h$ (the only term
with non-zero degree, by Lemma 3.3(i)).\qed\enddemo\vskip 12pt

\subheading{\S 3.2. The second blow-up} The center of the second blow-up is a
$\P^1$-bundle $B_1$ over $B$
$$\CD
B_1 @>i_1>>  V_1\\
@Vp_1 VV  @VVV\\
B @>i>> \P^8
\endCD$$
so classes on $B_1$ are combinations of (the pull-backs of) $k,h$ and $c_1(
\Cal O_{B_1}(-1))$; we call this latter $e$, and observe it is the pull-back
from $V_1$ of the class of the exceptional divisor $E_1$.
\proclaim{Lemma 3.5}

(i) ${p_1}_*e^i=\left\{\gathered
0\\ -1\\ -3k+2dh-6h\\ -6k^2+9dkh-27kh\\ 24dk^2h-72k^2h\endgathered\quad
\gathered i=0\\i=1\\i=2\\i=3\\i=4\endgathered\right.$\vskip 3pt

(ii) $c(N_{B_1}V_1)=(1+e)(1+k+dh-e)^3$\vskip 3pt

(iii) $i_1^*P_1=dk+dh-e$\endproclaim
\demo{Proof} (iii) is immediate, as $P$ is non-singular and pulls-back on $B$
to $dk+dh$ (Lemma 3.3 (iii)).

For (i) and (ii) we need to produce $B_1\subset E_1$ more explicitly as the
projectivization of a rank-2 subbundle of $N_B\P^8$.

First define for any $p\in \P^2$ a rank-8 subbundle $H_p$ of the trivial
bundle $B\times\C^9$ over $B$: if $F$ is a polynomial defining $C$, and
$(k,q)\in B$, $A\in\C^9=\text{Hom}(\C^3,\C^3)$, say
$$((k,q),A)\in H_p \iff \sum_{i=0}^2\left(\frac{\partial F}{\partial x_i}
\right)_qA(p)_i=0$$
where $A(p)_i$ is the $i$-th coordinate of $A(p)$. So the fiber of $H_p$ over
$q$ is the hyperplane of matrices $A\in\C^9$ such that $A(p)\in$ line tangent
to $C$ at $q$. Notice that the above equation has degree $d-1$ in the
coordinates of $q$: thus (denoting by $\C^9$ the trivial bundle $B\times\C^9$,
for short)
$$c_1\left(\frac{\C^9}{H_p}\right)=(d-1)h\quad.$$
Now restrict the Euler sequence for $\P^8$ to $B$ via $B\overset i\to
\hookrightarrow\P^8$: $H_p\subset \C^9$ determines a subbundle $\Cal H_p$ of
$i^*T\P^8$ and we have the following diagram of bundles over $B$ (suppressing
pull-back as usual)
$$\CD
   @.     0       @.        0         @.         0      @.\\
@.      @VVV              @VVV                @VVV         @.\\
0 @>>> \Cal O @>>>  H_p\otimes \Cal O_{\P^8}(1) @>>> \Cal H_p @>>> 0\\
@.       @|               @VVV                @VVV         @.\\
0 @>>> \Cal O @>>> \C^9\otimes \Cal O_{\P^8}(1) @>>>   T\P^8  @>>> 0\\
@.      @VVV              @VVV                @VVV         @.\\
   @.    0    @>>> \dfrac{\C^9}{H_p}\otimes \Cal O_{\P^8}(1) @=
\dfrac{T\P^8}{\Cal H_p}   @>>>  0\\
@.       @.               @VVV                @VVV         @.\\
   @.             @.       0          @.         0      @.
\endCD$$
from which it follows
$$c\left(\frac{T\P^8}{\Cal H_p}\right)=c\left(\frac{\C^9}{H_p}\otimes \Cal
O_{\P^8}(1)\right)=1+k+dh\quad.$$
Also, observe that each $\Cal H_p$ contains $TB$.

Now let $p_1,p_2,p_3$ be non-collinear points. A matrix has image contained in
a line if and only if it sends three non-collinear points to that line, thus
the intersection $H_{p_1}\cap H_{p_2}\cap H_{p_3}$ is the rank-6 bundle over
$B=\cP^2\times C$ whose fiber over $(k,q)\in B$ consists of all matrices whose
image is contained in the line tangent to $C$ at $q$. This is the space we
used to define $B_1$: if we set $\Cal Q= \Cal H_{p_1}\cap\Cal H_{p_2}\cap
\Cal H_{p_3}$, then
$$B_1=\P\left( \frac{\Cal Q}{TB}\right) \subset \P(N_B\P^8)=E_1,
\quad\text{and}\quad c\left(\frac{T\P^8}{\Cal Q}\right)= (1+k+dh)^3\quad.$$
Finally, the Euler sequences for $E_1$ and $B_1$ give the diagram
$$\CD
   @.     0      @.         0       @.           0      @.\\
@.      @VVV              @VVV                @VVV         @.\\
0 @>>> \Cal O @>>> \dfrac{\Cal Q}{TB}\otimes \Cal O_{B_1}(1) @>>> TB_1|B @>>>
0\\
@.       @|               @VVV                @VVV         @.\\
0 @>>> \Cal O @>>> N_B\P^8\otimes \Cal O_{B_1}(1) @>>>   TE_1|B  @>>> 0\\
@.      @VVV              @VVV                @VVV         @.\\
   @. 0 @>>> \dfrac{T\P^8}{\Cal Q}\otimes \Cal O_{B_1}(1) @=N_{B_1}E_1 @>>>
0\\
@.       @.               @VVV                @VVV         @.\\
   @.            @.        0        @.           0      @.
\endCD$$
(here $TB_1|B$, $TE_1|B$ denote the relative tangent bundles of $B_1$,
$E_1$ over $B$) from which
$$c(N_{B_1}E_1)=c\left(\frac{T\P^8}{\Cal Q}\otimes \Cal O_{B_1}(1)\right)
=(1+k+dh-e)^3\quad.$$
{From} this discussion, it's easy to obtain (i) and (ii):

$$\align
\text{(i) }{p_1}_*\sum_i&(-1)^ie^i=c\left(\frac{\Cal Q}{TB}\right)^{-1}\qquad
\qquad\qquad\text{by \cite{Fulton}, Proposition 3.1 (a)}\\
&=c\left(\frac{T\P^8}{\Cal Q}\right)c(N_B\P^8)^{-1}\qquad\qquad\quad
\text{by Whitney's formula}\\
&=\frac{(1+k+dh)^3(1+k)^3(1+h)^3}{(1+k+h)^9(1+dh)}\quad
\text{by the above and Lemma 3.3 (ii)}\\
&=1-3k+2dh-6h+6k^2-9dkh+27kh+24dk^2h-72k^2h\quad.
\endalign$$

(ii) $c(N_{B_1}V_1)=c(N_{E_1}V_1)c(N_{B_1}E_1)=(1+e)(1+k+dh-e)^3\quad.\qed$
\enddemo

Lemma 3.5 allows us to compute the term needed to apply Proposition 3.2 to the
second blow-up:
\proclaim{Lemma 3.6}
$$\int_{B_1}\frac{(B_1+i_1^*P_1)^8}{c(N_{B_1}V_1)}=
d(2d-3)(322d^2-1257d+1233)$$\endproclaim
\demo{Proof} This is
$$\int_{B_1}\frac{(1+dk+dh-e)^8}{(1+e)(1+k+dh-e)^3}$$
by Lemma 3.5 (ii) and (iii). Since the degree doesn't change after
push-forwards, this is also
$$\int_B{p_1}_*\frac{(1+dk+dh-e)^8}{(1+e)(1+k+dh-e)^3}\quad.$$
Computing the degree-4 term in the expansion of the fraction and applying
Lemma 3.5 (i) and the projection formula, this is computed as a sum of
degree-3 terms in $k,h$ over $B$. Lemma 3.3 (i) is used then to obtain the
stated expression.\qed\enddemo\vskip 12pt

\subheading{\S 3.3. The third blow-up} At this point we have to start taking
flexes into account. For any $q\in C$, let $f\ell(q)$ be the order of $q$ as a
flex of $C$, in the sense of \S 2.2: so $f\ell(q)=0$ if $q$ is {\sl not\/} a
flex of $C$, $f\ell(q)=1$ if $q$ is a {\sl simple\/} flex of $C$, and so
on.

The center $B_2 \overset i_2 \to\hookrightarrow V_2$ of the third
blow-up is the disjoint union
$$B_2=\bigcup_{f\ell(q)>0}\bq_2\quad,$$
where each $\bq_2$ maps isomorphically to the restriction $\bq_1$ of
the $\P^1$-bundle $B_1$ to $\cP^2\times\{q\}\subset B$. Moreover,
$B_2\cap \Til E_1 =\emptyset$ (Lemma 2.4). As $h$ restricts to 0 on
each $\cP^2\times \{q\}$, the intersection ring of $\bq_2$ is
generated by $k,e$ (defined as in \S 3.2).
Also, we denote by $e'$ the pull-back of $E_2$ to $\bq_2$, and by $p_{20}$ the
map $\bq_2 @>>> \cP^2\times\{q\}\cong \P^2$.
\proclaim{Lemma 3.7}

(i) $e'=e$

(ii) ${p_{20}}_*e^i=\left\{\gathered
0\\ -1\\ -3k\\ -6k^2\endgathered\quad\gathered
i=0\\ i=1\\ i=2\\ i=3\endgathered\right.$

(iii) $c(N_{\bq_2}V_2)=(1+e)(1+k-2e)^3$

(iv) $i_2^*P_2=dk-2e$\endproclaim

\demo{Proof} (ii) follows from Lemma 3.5 (i), since the restriction of $h$ to
$\bq_2$ is 0.

The key observation for the other points is that $\bq_2\cap \Til E_1=
\emptyset$. Realize $\bq_2\subset \P(N_{B_1}V_1)$ as $\P(\Cal L)$, where
$\Cal L$ is a sub-line bundle of $N_{B_1}V_1$. $\Til E_1\cap E_2$ is the
exceptional divisor of the blow-up of $E_1$ along $B_1$, i\.e\. the
projectivization of $N_{B_1}E_1$ in $N_{B_1}V_1$. That $\P(\Cal L)$ and
$\P(N_{B_1}E_1)$ are disjoint says that $\Cal L\cap N_{B_1}E_1$ is the
zero-section of $N_{B_1}V_1$, and therefore
$$\Cal L\cong \frac{N_{B_1}V_1}{N_{B_1}E_1}=N_{E_1}V_1\quad\text{ as bundles
on $\bq_1$.}$$

(i) With the same notations, $\Cal L$ is tautologically the universal line
bundle over $\P(\Cal L)$; it must then equal the restriction to $\bq_2$ of the
universal line bundle $\Cal O_{E_2}(-1)\cong N_{E_2}V_2$. In other words
$$\Cal L\cong N_{E_2}V_2\quad\text{ as bundles on $\bq_2$.}$$
Since the projection from $\bq_2$ to $\bq_1$ is an isomorphism, it follows
that
$$e=c_1(N_{E_1}V_1)=c_1(\Cal L)=c_1(N_{E_2}V_2)=e'\quad.$$

(iii) Call $\Eq_2$ the restriction of $E_2=\P(N_{B_1}V_1)$ to $\bq_1$. We have
Euler sequences
$$\CD
0 @>>> \Cal O @>>> \Cal L\otimes \Cal O(1) @>>> T\bq_2|\bq_1 @>>> 0\\
@. @VVV @VVV @VVV @.\\
0 @>>> \Cal O @>>> N_{B_1}V_1\otimes \Cal O(1) @>>> T\Eq_2|\bq_1 @>>> 0
\endCD$$
and we just argued $\Cal L\cong\Cal O(-1)$: so
$$\align
c(N_{\bq_2}\Eq_2)&=c\left(\frac{N_{B_1}V_1}{\Cal L}\otimes\check{\Cal L}\right)
\qquad\text{(restricted to $\bq_2$)}\\
&=\frac{(1+e-e')(1+k-e-e')^3}{(1+e'-e')}\\
&=(1+k-2e)^3\quad\text{by (i)};
\endalign$$
next, since $N_{\bq_1}B_1$ is clearly trivial, we have $c(N_{\Eq_2}E_2)=1$; so
putting $N_{\bq_2}V_2$ together:
$$c(N_{\bq_2}V_2)=c(N_{E_2}V_2)c(N_{\Eq_2}E_2)c(N_{\bq_2}\Eq_2)=
(1+e)(1+k-2e)^3\quad,$$
as claimed.

(iv) Since $P_1$ is non-singular along $B_1$, $P_2$ restricts to $dk-e-e'
=dk-2e$ by (i).\qed\enddemo

We are ready for the term needed to apply Proposition 3.2 to the third
blow-up:
\proclaim{Lemma 3.8}
$$\int_{B_2}\frac{(B_2+i_2^*P_2)^8}{c(N_{B_2}V_2)}= \sum_{f\ell(q)>0}
(196d^2-960d+1125)$$\endproclaim

\demo{Proof} By Lemma 3.7 (iii) and (iv), this is
$$\sum_{f\ell(q)>0}\int_{\bq_2}\frac{(1+dk-2e)^8}{(1+e)(1+k-2e)^3}
=\sum_{f\ell(q)>0}\int_{\P^2}{p_{20}}_*\frac{(1+dk-2e)^8} {(1+e)(1+k-2e)^3}
$$
(pushing forward doesn't change degrees) and one concludes with the projection
formula and Lemma 3.7 (ii).\qed\enddemo\vskip 12pt

\subheading{\S 3.4. Further blow-ups} Further blow-ups are necessary if there
are points $q$ on $C$ with $f\ell(q)>1$. We first attack the initial step.

The center $B_3 \overset i_3\to\hookrightarrow V_3$ of the fourth blow-up is
the disjoint union
$$B_3=\bigcup_{f\ell(q)>1} \bq_3\quad,$$
where each $\bq_3$ is a $\P^1$-bundle over $\bq_2$. The intersection ring of
$\bq_3$ is generated by (the pull-back of) the classes $k,e$ of $\bq_2$, and
by the class of the universal line bundle, i\.e\. the pull back $f$ of $E_3$
from $V_3$. Denote by $p_3$ the projection $\bq_3 @>>> \bq_2$.
\proclaim{Lemma 3.9}

(i) ${p_3}_*f^i=\left\{\gathered
0\\-1\\-e\\-e^2\\-e^3\endgathered\quad\gathered
i=0\\i=1\\i=2\\i=3\\i=4\endgathered\right.$

(ii) $c(N_{\bq_3}V_3)=(1+f)(1+k-2e-f)^3$

(iii) $i_3^*P_3=dk-2e-f$\endproclaim
\demo{Proof} (iii) is clear, as $P_2$ is non-singular along $\bq_3$.

For the other items, we have to produce $\bq_3\subset
\Eq_3=\P(N_{\bq_2}V_2)$ explicitly as the projectivization of a rank-2
subbundle of $N_{\bq_2}V_2$. Recall that each fiber of $\bq_3$ is
spanned by two points corresponding respectively to (1) a direction
transversal to $E_2$, and (2) a direction in $E_2$, transversal to
the fiber of $E_2$. Since these two points are always distinct, $\bq_3=\P(
\Cal L_1\oplus \Cal L_2)$, where $\P\Cal L_1$, $\P\Cal L_2$ give the two
distinguished points on each fiber. Now, $\Cal L_1\cap N_{\bq_2}E_2$ is the
zero-section in $N_{\bq_2}V_2$ (the first direction is transversal to $E_2$);
so, with $\Cal L$ as in the proof of 3.7,
$$\Cal L_1\cong N_{E_2}V_2\cong \Cal L\quad.$$
Similarly, since the second direction is transversal to the {\sl fiber\/} of
$E_2$, whose normal bundle in $E_2$ is trivial, $\Cal L_2\cong \Cal O$; and
therefore we have
$$\bq_3=\P(\Cal L\oplus\Cal O)\quad.$$

(i) As in the proof of 3.5 (i),
$${p_3}_*\sum_i(-1)^if^i=c(\Cal L\oplus\Cal O)^{-1}=\sum_i(-1)^ie^i$$
and (i) follows by matching dimensions.

(ii) Another pair of Euler sequences: on $\bq_3$
$$\CD
0 @>>> \Cal O @>>> (\Cal L\oplus\Cal O)\otimes \Cal O(1) @>>> T\bq_3|\bq_2
@>>> 0\\
@. @VVV @VVV @VVV @.\\
0 @>>> \Cal O @>>> N_{\bq_2}V_2\otimes \Cal O(1) @>>> T\Eq_3|\bq_2 @>>> 0
\endCD$$
Since $c_1(\Cal O(1))=-f$ and $E_3$ is the disjoint union of the $\Eq_3$:
$$\align
c(N_{\bq_3}E_3)&= c(N_{\bq_3}\Eq_3)\\
&=c\left(\frac{N_{\bq_2}V_2}{\Cal L\oplus\Cal O}\otimes \Cal O(1)\right)\\
&=(1+k-2e-f)^3
\endalign$$
(the Chern roots of $N_{\bq_2}V_2$ are $e,k-2e,k-2e,k-2e,0$ by Lemma 3.7
(iii)). Finally:
$$c(N_{\bq_3}V_3)=c(N_{E_3}V_3)c(N_{\bq_3}E_3)=(1+f)(1+k-2e-f)^3$$
as stated.\qed\enddemo

Lemma 3.9 describes the situation at the fourth blow-up. The next
blow-ups are built on this in the sequence described in \S 2.4: the
center $B_j \overset i_j\to \hookrightarrow V_j$ of the $(j+1)$-st
blow-up ($j\ge 3$) is the disjoint union
$$B_j=\bigcup_{f\ell(q)> j-2} \bq_j\quad,$$
where each $\bq_j$ maps isomorphically down to $\bq_3$, and is
disjoint from $\Til E_{i-1}$ (Lemma 2.11). The intersection ring of
each $\bq_j\cong\bq_3$ is then generated by $k,e,f$, and the relations
stated in Lemma 3.9 (i) hold, for the projection $p_{j2}:
\bq_j @>>> \bq_2$. Denote by $f_j$ the pull-back of $E_j$ to $\bq_j$;
Lemma 3.9 can be extended to all stages in the sequence:

\proclaim{Lemma 3.9 (continued)} For $3\le j\le f\ell(q)+1$

$(i)_j$ $f_j=f$

$(ii)_j$ $c(N_{\bq_j}V_j)=(1+f)(1+k-2e-(j-2)f)^3$

$(iii)_j$ $i_j^*P_j=dk-2e-(j-2)f$\endproclaim

\demo{Proof} For $j=3$ this is given by Lemma 3.9. So it suffices to show
that, for $3\le j\le f\ell(q)$, $(i)_j,(ii)_j,(iii)_j$ imply $(i)_{j+1},
(ii)_{j+1}, (iii)_{j+1}$. Consider then $\bq_{j+1}=\P(\Cal L_{j+1}) \subset
\P(N_{\bq_j}V_j)$. So $f_{j+1}$ is the class of $\Cal O_{\bq_{j+1}}(-1)$,
i\.e\. of $\Cal L_{j+1}$. Since $\bq_{j+1}\cap\Til E_j=\emptyset$ (Lemma 2.11
(iv)), we get by the usual argument
$$f_{j+1}=c_1(\Cal L_{j+1})=c_1(N_{E_j}V_j)=f_j\quad:$$
and $f_j=f$ by $(i)_j$; so $f_{j+1}=f$, giving $(i)_{j+1}$.

$(iii)_{j+1}$ follows then from $(iii)_j$ and $(i)_{j+1}$, since $P_j$ is
non-singular along $B_j$.

Finally, we use the Euler sequences
$$\CD
0 @>>> \Cal O @>>> \Cal L_{j+1}\otimes \Cal O(1) @>>> T\bq_{j+1}|\bq_j
@>>> 0\\
@. @VVV @VVV @VVV @.\\
0 @>>> \Cal O @>>> N_{\bq_j} V_j\otimes \Cal O(1) @>>> T\Eq_{j+1}|\bq_j
@>>> 0
\endCD$$
to get (since $E_{j+1}$ is the disjoint union of the $\Eq_{j+1}$)
$$\align
c(N_{\bq_{j+1}}E_{j+1})&=c(N_{\bq_{j+1}}\Eq_{j+1})\\
&=c\left(\frac{N_{\bq_j} V_j}{\Cal L_{j+1}}\otimes\Cal O(1)\right)\\
&=\frac{(1+f-f)(1+k-2e-(j-2)f-f)^3}{(1+f-f)}\qquad\text{by $(ii)_j$}\\
&=(1+k-2e-(j-1)f)^3\qquad,
\endalign$$
so
$$c(N_{\bq_{j+1}}V_{j+1})=c(N_{E_{j+1}}V_{j+1})c(N_{\bq_{j+1}}E_{j+1})=
(1+f)(1+k-2e-(j-1)f)^3\quad,$$
i\.e\. $(ii)_{j+1}$.\qed\enddemo

We get then the key term to apply Proposition 3.2 to the $j$-th blow up in the
sequence. In fact, we can cover Lemma 3.8 as well in one statement:
\proclaim{Lemma 3.10} For $j\ge 2$
$$\multline\int_{B_j}\frac{(B_j+i_j^*P_j)^8}{c(N_{B_j}V_j)}=\sum_{f\ell(q)>j-2}
30j^4-96(d-1)j^3\\
+12(d-1)(7d-11)j^2+84(d-1)^2j-7(2d-3)(22d-39).\endmultline$$
\endproclaim
\demo{Proof} For $j=2$, this is Lemma 3.8. For $j\ge 3$, by Lemma 3.9 this is
$$\sum_{f\ell(q)>j-2}\int_{\bq_j}\frac{(1+dk-2e-(j-2)f)^8}
{(1+f)(1+k-2e-(j-2)f)^3}\quad.$$
If $p_{j2}$ denotes the projection $\bq_j @>>> \bq_2$, (and $p_{20}$ is the
map $\bq_2 @>>> \cP^2\times\{q\}\cong\P^2$, as in \S 3.3), this can be
computed as
$$\sum_{f\ell(q)>j-2}\int_{\P^2}{p_{20}}_*{p_{j2}}_*\frac
{(1+dk-2e-(j-2)f)^8}{(1+f)(1+k-2e-(j-2)f)^3}\quad,$$
which is evaluated by using the projection formula, 3.9 (i) and 3.7 (ii).
\qed\enddemo\vskip 12pt

\subheading{\S3.5. The predegree of $\ocbar$} Computing the predegree
of $\ocbar$ is now a straightforward application of Proposition 3.2
and Lemmas 3.4, 3.6 and 3.10: by Proposition 3.2
$$\Til P^8=P^8-\sum_{j\ge 0}\int_{B_j}\frac{(B_j-i_j^*P_j)^8}{c(N_{B_j}V_j)}$$
(where $B_0=B$, etc\.), and the terms in the summation have been computed in
sections 3.1--3.4. This gives

\proclaim{Proposition 3.11} The predegree of $\oc$ is
$$\multline
d^8-d(10d-9)(14d^2-33d+21)-d(2d-3)(322d^2-1257d+1233)\\
-\sum_{j\ge 2}\sum
\Sb q\in C\\ f\ell(q)>j-2\endSb
(30j^4-96(d-1)j^3+12(d-1)(7d-11)j^2\\
+84(d-1)^2j-7(2d-3)(22d-39)).\endmultline$$\endproclaim

This result can be given in handier forms. For example:
\proclaim{Theorem III(a)} The predegree of $\oc$ is
$$\gather
d(d-2)(d^6+2d^5+4d^4+8d^3-1356d^2+5280d-5319)-\sum_{q\in C}f\ell(q)
(f\ell(q)-1)\\
\left(6f\ell(q)^3+(75-24d)f\ell(q)^2+(28d^2-240d+393)f\ell(q)+196d^2-
960d+1125\right)
\endgather$$\endproclaim

\demo{Proof} Invert the order of the summations in Proposition 3.11,
then
use the fact that $\sum_{q\in C} f\ell(q)=3d(d-2)$ (the number of
flexes of $C$, counted with multiplicity).
\qed\enddemo

Or, in another form:
\proclaim{Theorem III(b)} Denote by $f_C^{(r)}$ the sum
$\sum_{q\in C}f\ell(q)^r$. Then the predegree of $\oc$ is
$$\multline
d^8-8d(98d^3-492d^2+843d-486)-(168d^2-720d+732)f_C^{(2)}\\
-(28d^2-216d+318)f_C^{(3)}-(69-24d)f_C^{(4)}-6f_C^{(5)}.\endmultline
$$\endproclaim\vskip 12pt

By Theorem III(B), if $C$ is smooth then {\sl the predegree of $\oc$
depends only on the degree $d$ of $C$ and on the four numbers
$f_C^{(2)}$, $f_C^{(3)}$, $f_C^{(4)}$ and $f_C^{(5)}$.\/}\vskip 12pt

If $C$ only has simple flexes, then $f\ell(q)=0$ or 1 for all $q\in C$, so
Theorem III(A) gives

\proclaim{Corollary} If all flexes of $C$ are simple, then the
predegree of $\oc$ is
$$\multline d(d-2)(d^6+2d^5+4d^4+8d^3-1356d^2+5280d-5319)\\
=d^8-1372d^4+7992d^3-15879d^2+10638d\quad.\endmultline$$
\endproclaim
Denoting this polynomial in $d$ by $P(d)$, we remark that it gives
the degree of the orbit closure of the general smooth
plane curve of degree $d\ge 4$ (indeed, such a curve $C$ has no non-trivial
automorphisms, so by Lemma 3.1 the degree of $\ocbar$ equals the predegree).

{\smc Remark.} Denoting by $f_k(d)$ the (negative) contribution to the
predegree arising from a flex
of order $k$ on a curve of degree $d$, we have, as an immediate
consequence of Theorem III(A):
$$
f_k(d)=-k(k-1)((28k+196)d^2-(24k^2+240k+960)d+(6k^3+75k^2+393k+1125)).
$$
E.g., $f_2(d)=-6(84d^2-512d+753)$ and $f_3(d)=-6(280d^2-1896d+3141)$.
It is an easy calculus exercise to show that $f_k(d)<0$ whenever $d\ge
k+2\ge4$. This proves that the predegree is maximal for a curve with
only simple flexes.\vskip 12pt

\subheading{\S3.6. Examples} It is a consequence of Lemma 3.1 that the
predegree of the orbit of a smooth plane curve is divisible by the order of
its $\PGL(3)$@-stabilizer. This cuts both ways. On the one hand, each curve
with non-trivial automorphisms provides us with a non-trivial check of the
formulas above. On the other hand, these formulas might help in determining
which automorphism groups of smooth plane curves occur. We illustrate this
below.

Consider, for $d\ge3$, the Fermat curve $x^d+y^d+z^d$.
Its $3d$ flexes have order $d-2$, so
the predegree of its orbit is $P(d)+3d\cdot f_{d-2}(d)$. So for each $d$
this {\sl number\/} is divisible by $6d^2$,
the order of the stabilizer. This implies that in the ring
$\Z[d]$ the {\sl polynomial\/} $P(d)+3d\cdot f_{d-2}(d)$ is divisible by
$d^2$ and
that the quotient polynomial takes values divisible by 6. Indeed
$$P(d)+3d\cdot f_{d-2}(d) = d^2(d-2)(d^5+2d^4-26d^3-7d^2+192d-192)\quad.$$
Dividing this by $6d^2$, we get the degree of the orbit closure of the Fermat
curve, i.e., of the trisecant variety to the $d$@-uple embedding of
$\P^2$ in $\P^N$, as mentioned in the introduction.

Here is a similar example for all $d\ge5$: the curve
$x^{d-1}y+y^{d-1}z+z^{d-1}x$.
The points $(1:0:0)$, $(0:1:0)$ and $(0:0:1)$ are flexes of order $d-3$;
counted
with multiplicity, $3(d^2-3d+3)$ flexes remain. The group $D$ of diagonal
matrices with entries $(1,\zeta,\zeta^{2-d})$, where $\zeta$ is a
$(d^2-3d+3)$@-rd root of unity, acts on the latter flexes without fixed
points; so either there is one orbit of flexes of order 3, or one orbit
of flexes of order 2 and one orbit of simple flexes, or, finally,
three orbits of simple flexes. Now one uses the automorphism
$\sigma\colon(x:y:z)\mapsto(y:z:x)$ to exclude the first two possibilities;
moreover, one verifies that the automorphism group $G$ of the curve is the
semidirect product of $D$ and $<\sigma>$ (and that the simple flexes form
one $G$@-orbit). So the degree of the orbit closure is
$$\frac{\vphantom{\int}P(d)+3f_{d-3}(d)}{\vphantom{\int}3(d^2-3d+3)}
=\frac13(d^6+3d^5+6d^4-21d^3-1354d^2+5463d-5508)\quad.$$

Next we list, for some small values of $d$, the
numbers (and their factorizations) we get from the corollary to Theorem
III:
\protable
&&$d$&&$\vphantom{\int}P(d)$&&$\vphantom{\int}P(d)$ factored&\cr
\tablerule
&&3 && 216&& $\vphantom{2^{\{}}2^3\cdot3^3$&\cr\tablerule
&&4 && 14280&& $\vphantom{2^{\{}}2^3\cdot3\cdot5\cdot7\cdot17$&\cr\tablerule
&&5 && 	188340&& $\vphantom{2^{\{}}2^2\cdot3\cdot5\cdot43\cdot73$&\cr\tablerule
&&6 && 1119960&&
$\vphantom{2^{\{}}2^3\cdot3^3\cdot5\cdot17\cdot61$&\cr\tablerule
&&7 && 4508280&&
$\vphantom{2^{\{}}2^3\cdot3^2\cdot5\cdot7\cdot1789$&\cr\tablerule
&&8 && 14318256&& $\vphantom{2^{\{}}2^4\cdot3\cdot317\cdot941$&\cr\tablerule
&&9 && 38680740&&
$\vphantom{2^{\{}}2^2\cdot3^6\cdot5\cdot7\cdot379$&\cr\tablerule
&&10 && 92790480&& $\vphantom{2^{\{}}2^4\cdot3\cdot5\cdot59\cdot6553$&\cr.

So for $d=3$ we get 216 for the predegree of the orbit of any smooth
plane cubic curve. This gives the well-known numbers 12, resp.~6,
resp.~4 for the degree of the orbit closure of a smooth plane cubic
with $j\ne0, 1728$, resp.~$j=1728$, resp.~$j=0$.
Note that the group of
projective automorphisms of a smooth cubic contains the 9 translations over
points of order dividing 3 as a normal subgroup. The quotient can be
identified
with the automorphisms that fix a given flex. Thus there exist 18, resp.~36,
resp.~54 projective automorphisms when $j\ne0, 1728$, resp.~$j=1728$,
resp.~$j=0$.

For $d=4$ we get 14280 for the predegree of the orbit of a smooth
plane quartic with only simple flexes. An example of such a curve is
the Klein curve $x^3y+y^3z+z^3x$; it has 168 automorphisms, so the
degree of its orbit closure is $14280/168=85$.

If a smooth quartic has $n$ hyperflexes (i.e., flexes of order 2), the
predegree of its orbit
equals $14280-294n$. E.g., the degree of the orbit closure of the
Fermat quartic is 112, as there are 12 hyperflexes and 96
automorphisms. As an other example, consider the curve
$x^4+xy^3+yz^3$. It has 1 hyperflex and 9 automorphisms, so the degree
of its orbit closure is $(14280-294)/9=1554$.  In fact, in
\cite{Vermeulen} there is a complete list of the automorphism groups
that occur for a quartic with a given number of hyperflexes. The
implied congruence conditions are equivalent to requiring that $P(4)$
be divisible by 168 and that $P(4)+28f_2(4)$ be divisible by 2016.
(This follows already from the existence of the 3 quartics above.)

In the other direction, these formulas give non-trivial information on
the automorphism groups of plane curves. Consider smooth plane curves
of degree $d$ with only simple flexes. The least common multiple of
the orders of the stabilizers of these curves divides $P(d)$. Now it
is well-known that a smooth curve of positive genus
$g(={{d-1}\choose2})$ cannot have an automorphism of prime order
$p>2g+1(=d^2-3d+3)$. Using the Hurwitz formula one also excludes the
cases $(d,g,p)=(4,3,5)$, $(6,10,17)$ and $(10,36,59)$. Looking at the
table above, we conclude then that said l.c.m. divides 216, 168, 60,
1080, 2520, 48, 102060, 240 respectively for $d$ equal to 3, 4, 5, 6,
7, 8, 9, 10 respectively.

These bounds seem to be pretty good: by the above, the actual l.c.m.
equals 108 (resp.~168) for $d=3$ (resp.~4); it's not unreasonable to
expect that the bound is sharp for $d=5$, 8 and 10 (perhaps there even
exist curves with automorphism groups of this order); the Valentiner
sextic has only simple flexes and 360 automorphisms (cf. \cite{BHH}),
so the bound for $d=6$ is sharp if and only if there exists a sextic
with only simple flexes and with 27 dividing the order of its
stabilizer. Finally, for $d=9$ the bound is probably not optimal.

\enddocument